\documentclass[11pt]{article}

\usepackage[a4paper,margin=2.7cm]{geometry}
\usepackage[T1]{fontenc}
\usepackage[utf8]{inputenc}

\usepackage{amsmath,amssymb,amsfonts,amsthm}
\usepackage{mathtools}
\usepackage{mathrsfs}
\usepackage{bm}
\usepackage{bbm}

\usepackage{graphicx}
\usepackage{xcolor}
\usepackage{hyperref}

\usepackage{tikz}
\usetikzlibrary{arrows.meta,positioning,shapes.geometric}

\theoremstyle{definition}
\newtheorem{defi}{Definition}[section]

\theoremstyle{plain}
\newtheorem{prop}[defi]{Proposition}

\newtheorem{theorem}[defi]{Theorem}

\theoremstyle{remark}
\newtheorem{rem}[defi]{Remark}

\newcommand{\ii}{\mathrm{i}}

\newcommand{\CC}{\mathbb{C}}
\newcommand{\HH}{\mathbb{H}}
\newcommand{\OO}{\mathbb{O}}

\newcommand{\buTwo}{\mathbbm{1}_2}
\newcommand{\rev}{\boldsymbol{\curvearrowleft}}

\newcommand{\ent}{\mathbin{\rotatebox[origin=c]{90}{\scriptsize$\boldsymbol{\circlearrowleft}$}}}

\newcommand{\st}{^{\!*}}

\title{A Mnemonic Matrix Rule for (Split) Octonionic Multiplication\\
and its Extension to the Cayley--Dickson Tower}
\author{Jean-Pierre Gazeau\\
Universit\'e Paris Cit\'e, CNRS, Astroparticule et Cosmologie, F-75013 Paris, France\\
Faculty of Mathematics, University of Bia\l ystok, 15-245 Bia\l ystok, Poland\\
\small\textit{gazeau@apc.in2p3.fr, j.gazeau@uwb.edu.pl}}
\date{\today}

\begin{document}
\maketitle

\begin{abstract}
We present a compact mnemonic device for computing the product of two
(split) octonions written in Cayley--Dickson form $q+\ell p$ with $q,p\in\HH$.
The rule appears as a simple $(R+L)$ pattern of right-ordered and left-ordered
(quaternionic) products inside a $2\times2$ quaternionic matrix model.
The pattern extends verbatim to all algebras in the Cayley--Dickson tower,
providing an efficient computational tool in non-associative settings.
To our knowledge, this explicit ``$(R+L)$'' mnemonic does not appear in the
classical literature on octonions or composition algebras.
\end{abstract}

\section{Introduction}

The octonions $\OO$ \cite{Baez,Schafer,Okubo,SpringerVeldkamp,ConwaySmith} may be described as pairs of quaternions
$(q,p)\in \HH\oplus \HH\,\ell$ via the Cayley--Dickson construction.
The multiplication is
\begin{equation}
(q_1+\ell p_1)(q_2+\ell p_2)
= q_1 q_2 + \ell^2\, p_2 p_1\st \;+\; \ell\,(q_2 p_1 + q_1\st p_2),
\label{eq:CD}
\end{equation}
where $\st$ denotes quaternionic conjugation and $\ell^2=\pm1$ covers the
division and split cases.
While \eqref{eq:CD} is standard, its direct use can become cumbersome because it
mixes conjugation with order-sensitive products.

In this note we isolate a simple mnemonic rule that organizes octonion
multiplication into a $2\times2$ block pattern of quaternionic products, where
exactly half of the elementary products appear in reversed order.
The pattern propagates inductively to all Cayley--Dickson algebras
$\mathbb{A}_n$ of dimension $2^n$ and motivates a compatible ``interlaced''
product for $2\times2$ matrices over a unital non-associative algebra.

\section{Cayley--Dickson form and a $2\times2$ quaternionic model}

Write a (split) octonion as
\[
o = q+\ell p,\qquad q,p\in\HH,\qquad \ell^2=\pm1.
\]
Define the linear map
\begin{equation}
M(o):=
\begin{pmatrix}
 q & \ell^2\, p\st\\
 p & q\st
\end{pmatrix}.
\label{eq:M}
\end{equation}

\begin{defi}[Reversal marker]
Given a product $xy$ in a (possibly non-commutative) algebra, we denote by
$\overset{\rev}{xy}$ the reversed product $yx$.
\end{defi}
\begin{rem}
The symbol $\rev$ is purely typographical: it does not define a new operation
on the algebra. It only indicates that, in the marked contraction, the two
factors are multiplied in the opposite order.
\end{rem}

\begin{defi}[Interlaced $2\times2$ product over $\HH$]
For matrices of the form \eqref{eq:M}, define an ``interlaced'' multiplication
$\ent$ by expanding the usual matrix product, but reversing the order in the
\emph{off-diagonal} contractions:
\begin{equation}
\begin{pmatrix}
 a & b\\ c & d
\end{pmatrix}
\ent
\begin{pmatrix}
 a' & b'\\ c' & d'
\end{pmatrix}
:=
\begin{pmatrix}
 aa' + \overset{\rev}{bc'} & \overset{\rev}{ab'} + bd'\\
 \overset{\rev}{ca'} + dc' & cb' + \overset{\rev}{dd'}
\end{pmatrix}.
\label{eq:ent-general}
\end{equation}
\end{defi}

A direct computation shows that the Cayley--Dickson multiplication is reproduced
at the level of $M(\cdot)$.

\begin{prop}[Matrix mnemonic for octonion multiplication]
For $o_i=q_i+\ell p_i$ with $q_i,p_i\in\HH$,
\begin{equation}
M(o_1)\ent M(o_2)=M(o_1o_2),
\label{eq:mnemo}
\end{equation}
where $o_1o_2$ is given by \eqref{eq:CD}.
\end{prop}

\begin{rem}
The point of \eqref{eq:mnemo} is not a new model of $\OO$ (many exist), but a
\emph{visible order-pattern}: each block in the product is a sum of a
right-ordered product and a left-ordered one.
\end{rem}
\begin{rem}
Unlike Zorn's vector--matrix formulation \cite{ConwaySmith}, where non-associativity is encoded in
the replacement of ordinary products by vector dot and cross products in the
off--diagonal blocks, the present construction retains standard block entries
and shifts all non-associativity into a systematic pattern of order reversals
within an interlaced $2\times2$ matrix multiplication.
\end{rem}

\section{The mnemonic $(R+L)$ pattern}

Inspection of \eqref{eq:ent-general}--\eqref{eq:mnemo} yields the rule:

\begin{itemize}
\item each block is the sum of exactly two products;
\item one term is in \emph{right} (natural) order (call it $R$);
\item the other is in \emph{left} (reversed) order (call it $L$);
\item the pattern is
\[
\begin{pmatrix}
 R+L & L+R\\
 L+R & R+L
\end{pmatrix}.
\]
\end{itemize}

This compact $(R+L)$ pattern provides a fast and reliable way to compute
octonionic products in Cayley--Dickson form without repeatedly re-deriving
where conjugations and order-reversals must occur.

\section{Propagation through the Cayley--Dickson tower}

Let $\mathbb{A}_n$ denote the $2^n$-dimensional algebra obtained by iterating the
Cayley--Dickson doubling.
Assume $\mathbb{A}_n$ is equipped with its standard involution $\st$.
Represent an element of $\mathbb{A}_{n+1}$ as $a+\epsilon b$ with
$a,b\in\mathbb{A}_n$ and $\epsilon^2=\pm1$, and define
\[
M_{n+1}(a+\epsilon b):=
\begin{pmatrix}
 a & \epsilon^2\, b\st\\
 b & a\st
\end{pmatrix}.
\]

\begin{theorem}[$(R+L)$ propagation]
The interlaced product \eqref{eq:ent-general} applied to $M_{n+1}(\cdot)$
reproduces the Cayley--Dickson multiplication on $\mathbb{A}_{n+1}$.
Consequently the $(R+L)$ order-reversal pattern propagates inductively through
all Cayley--Dickson algebras $\mathbb{A}_n$.
\end{theorem}

\begin{rem}
For $n\ge4$ the algebras cease to be alternative, and zero divisors appear.
The mnemonic remains valid as an \emph{identity in the defining multiplication},
but invertibility and ``determinant-like'' constructions must be stated with
care.
\end{rem}

\section{Toward matrix structures over non-associative algebras}

The representation \eqref{eq:M} may be viewed as the natural octonionic analogue
of the classical $2\times2$ real and complex matrix models for $\CC$ and $\HH$:
\begin{equation}
\label{eq:CHO}
\CC\ni z=x+\ii y\;\longmapsto\;
\begin{pmatrix} x & -y\\ y & x\end{pmatrix},
\qquad
\HH\ni q=q_0+q_1\mathbf{i}+q_2\mathbf{j}+q_3\mathbf{k}
\;\longmapsto\;
\begin{pmatrix}
 q_0+\ii q_3 & -q_2+\ii q_1\\
 q_2+\ii q_1 &  q_0-\ii q_3
\end{pmatrix}.
\end{equation}

At the level of the tower, one obtains maps
\[
M_n:\mathbb{A}_n\longrightarrow \mathrm{Mat}_2(\mathbb{A}_{n-1}),
\]
and the interlaced product provides a natural way to retain a controlled
multiplication on these $2\times2$ blocks even when $\mathbb{A}_n$ is
non-associative. In particular, one may introduce:

\begin{itemize}
\item \emph{flexible or alternative matrix products}, obtained by restricting
the matrix multiplication on $\mathrm{Mat}_2(\mathbb{A}_{n-1})$ to combinations
that respect the underlying (left, right, or fully) alternative laws in
$\mathbb{A}_n$;

\item \emph{Moufang-type identities on special subspaces}, such as the
``Hermitian'' or ``pure'' components of $\mathbb{A}_n$, where the algebra
retains sufficient alternativity for the classical Moufang identities to
descend to matrix form;

\item \emph{groups of norm-preserving transformations}, defined as those
$2\times 2$ matrices over $\mathbb{A}_{n-1}$ whose action (by left or right
multiplication, or by conjugation) preserves the quadratic form induced by the
Cayley--Dickson norm on $\mathbb{A}_n$;

\item \emph{generalized determinant-like invariants}, obtained by replacing the
usual determinant with quadratic or quartic expressions that remain invariant
under the above norm-preserving groups (the determinant of a Hermitian
$3\times 3$ octonionic matrix being the archetypal example).
\end{itemize}

Such constructions arise naturally in exceptional geometry and the theory of
Jordan algebras (Albert algebra), in the triality phenomena associated with
$\mathrm{Spin}(8)$, and in Tits' constructions of exceptional Lie groups.
Within this framework, the representation $M_n$ serves as a bridge between the
non-associative algebraic world of the Cayley--Dickson sequence and the
classical matrix formalisms underlying these exceptional structures.

\section{Interlaced $2\times2$ multiplication over a unital non-associative algebra}

Let $\mathcal{A}$ be a \emph{unital} (not necessarily associative, not necessarily
commutative) algebra. Define $\ent$ on $\mathrm{Mat}_2(\mathcal{A})$ by the same
rule \eqref{eq:ent-general}:
\begin{equation}
\label{eq:mnemo-block2}
\begin{pmatrix}
  a & b\\ c & d
\end{pmatrix}
\ent
\begin{pmatrix}
  a' & b'\\ c' & d'
\end{pmatrix}
:= \begin{pmatrix}
  aa^{\prime} +   \overset{\rev}{bc^{\prime}}  &  \overset{\rev}{ab^{\prime}}+ bd^{\prime}  \\
   \overset{\rev}{ca^{\prime}}  + dc^{\prime} &  cb^{\prime}+\overset{\rev}{dd^{\prime}}
\end{pmatrix}=
\begin{pmatrix}
  aa' + c'b & b'a + bd'\\
  a'c + dc' & cb' + d'd
\end{pmatrix},
\end{equation}
where the displayed form uses $\overset{\rev}{xy}=yx$ to remove the reversal
notation.

\begin{rem}[On inverses and ``determinants'']
In a general non-associative algebra $\mathcal{A}$, expressions such as
$ad-cb$ and $da-cb$ need not coincide, need not be central, and may fail to
associate with the surrounding factors. Therefore one should distinguish
\emph{left} and \emph{right} inverse candidates, and state sufficient
conditions ensuring that scalar prefactors can be moved safely.
\end{rem}

A convenient safe statement is the following.

\begin{prop}[One-sided inverses under nucleus/centrality assumptions]
Let $X=\begin{pmatrix}a&b\\ c&d\end{pmatrix}\in\mathrm{Mat}_2(\mathcal{A})$ and
define
\[
\Delta_L:=da-cb,\qquad \Delta_R:=ad-cb,\qquad
X^\sharp:=\begin{pmatrix} d & -b\\ -c & a\end{pmatrix}\,.
\]
Assume:
\begin{itemize}
\item $\Delta_L$ is invertible and belongs to the \emph{left nucleus} of
$\mathcal{A}$ (so that $(\Delta_L^{-1}x)y=\Delta_L^{-1}(xy)$ for all $x,y$),
\item $\Delta_R$ is invertible and belongs to the \emph{right nucleus} of
$\mathcal{A}$,
\end{itemize}
then
\[
(\Delta_L)^{-1} X^\sharp \ent X = \buTwo\,,
\qquad
X \ent X^\sharp (\Delta_R)^{-1} = \buTwo\,
\]
i.e.\ the displayed matrices provide a left inverse and a right inverse of $X$
for the interlaced product, yielding a loop structure \cite{Pflugfelder}.
\end{prop}

\begin{rem}
For alternative division algebras (e.g.\ $\HH$ and $\OO$), nucleus conditions are
automatically much easier to satisfy on many natural subalgebras, and one
recovers familiar inverse formulae in Zorn-type models \cite{daboul99,ConwaySmith}.
\end{rem}

\section{Conclusion}

The $(R+L)$ mnemonic yields a compact and effective method for computing with
octonions (and split octonions) in Cayley--Dickson form, and it propagates
uniformly through the entire Cayley--Dickson tower.
The same order-reversal pattern motivates an interlaced product on $2\times2$
matrices over a unital non-associative algebra, which in turn suggests
structured ``matrix-like'' formalisms adapted to non-associativity.
We hope this note will be useful both pedagogically and as a computational
device in exceptional algebra and geometry.

\subsection*{Acknowledgements} The author is grateful to Dr. Hamed Pehjan for his careful reading of the manuscript and for several helpful suggestions and corrections. 


\end{document}